\documentclass[number,citesort,dvips]{arxbj}

\aid{0}
\volume{17}
\issue{2}
\pubyear{2011}
\firstpage{628}
\lastpage{642}
\doi{10.3150/10-BEJ285}

\makeatletter
\newtheorem{proposition}{Proposition}
\newtheorem{corollary}{Corollary}

\def\bI{\mathbf{I}}
\def\E{\mathrm{E}}
\def\Var{\operatorname{Var}}
\def\tr{\operatorname{tr}}
\makeatother

\begin{document}
\begin{frontmatter}

\title{An approximate quantum Cram\'{e}r--Rao bound based on skew information}
\runtitle{Skew information and quantum Cram\'{e}r--Rao-type bounds}

\begin{aug}
\author{\fnms{Alessandra} \snm{Luati}\corref{}\ead[label=e1]{alessandra.luati@unibo.it}}

\runauthor{A. Luati}

\address{Department of Statistics,
University of Bologna,
via Belle Arti 41, 40126 Bologna, Italy.\\
\printead{e1}}
\end{aug}

\received{\smonth{7} \syear{2009}}
\revised{\smonth{4} \syear{2010}}

%
\begin{abstract}
A closed-form expression for Wigner--Yanase skew information in
mixed-state quantum systems is derived. It is shown that limit
values of the mixing coefficients exist such that Wigner--Yanase
information is equal to Helstrom information. The latter
constitutes an upper bound for the classical expected Fisher
information, hence the inverse Wigner--Yanase information provides
an approximate lower bound to the variance of an unbiased
estimator of the parameter of interest. The advantage of
approximating Helstrom's sharp bound lies in the fact that
Wigner--Yanase information is straightforward to compute, while it
is often very difficult to obtain a feasible expression for
Helstrom information. In fact, the latter requires the solution of
an implicit second order matrix differential equation, while the
former requires just scalar differentiation.
\end{abstract}

%
\begin{keyword}
\kwd{Cram\'{e}r--Rao-type bounds}
\kwd{Fisher information}
\kwd{parametric quantum models}
\end{keyword}

\end{frontmatter}

\section{Introduction}\label{sec:int}

Classical statistics applied to quantum systems gives rise to more
than one Fisher information quantity
for the unknown parameter $\theta$ that specifies the state of
the system, denoted by $\rho(\theta)$.
The classical Fisher information $i(\theta,M)$ measures the precision of an unbiased
estimator $t(x)$ of $\theta$ based on the outcome of an arbitrary measurement $M$, via
the Cram\'{e}r--Rao bound
\begin{equation}\label{eq:crb}
\Var\{t(x)\} \geq i(\theta,M) ^{-1}.
\end{equation}
Various quantum analogues of classical information may be obtained directly
based on quantum operators and without performing any measurement.
The duality between classical and quantum information
is a consequence of the interaction between the microscopic
environment, where quantum systems evolve, and the macroscopic
world, where measurements are performed. On the other hand, the
existence of many quantum versions of one classical quantity is a
characteristic feature of quantum mechanics, where states and
measurements are represented by non-commutative operators or, in
finite dimensions, by matrices.

Quantum analogs of classical information are derived according
to which expression of Fisher information is generalized to the
quantum setting, or to which operator-valued version of the
logarithmic derivative is chosen, for example, the left or
right logarithmic derivative (\cite{YuenLax1973}, \cite{Holevo1982}, Section 6.6) or the symmetric logarithmic derivative
(\cite{Helstrom1967,Helstrom1967}, Section 8.4). As a result, each quantum information
quantity inherits different properties and, consequently, may find
application in different inferential problems, essentially related
to Cram\'{e}r--Rao-type bounds. Examples are the Helstrom \cite{Helstrom1967}
bound, where the information is based on the symmetric logarithmic
derivative, and the two variants based on the left and right
logarithmic derivatives considered by Belavkin \cite{Belavkin1976}.
Holevo \cite{Holevo1982} derived a bound based on vectors of matrix derivatives that
can be attained asymptotically~\cite{GutaKahn2006}.

The relevance of the symmetric logarithmic derivative in the
statistical literature is due to its relation with classical
Fisher information. In fact, Helstrom information,
$I_H(\theta)$, constitutes an upper bound for classical Fisher
information \cite{braunsteincaves1994},
%
\begin{equation}i(\theta,M)\leq
I_H(\theta),\label{eq:infineq}
\end{equation}
from which the bound
%
\begin{equation}\Var\{t(x)\} \geq I_H(\theta) ^{-1}\label{eq:qrb}
\end{equation}
\cite{Helstrom1967} is obtained as a
corollary. Since, in the quantum setting, the classical
information formally depends on the measurement carried out in the
system, much attention has been devoted in the literature to the
search for measurements maximizing Fisher information. Major
contributions have come from \cite{BarndorffNielGill2000} and
\cite{Luati2004,Luati2008}, where necessary and sufficient conditions for
equality between classical Fisher and quantum Helstrom information
are derived, based on different hypotheses on the state of the
quantum system.

Another strand of the literature has investigated the relations
among quantum information quantities. Hayashi \cite{Hayashi2002} compares the
Kubo--Mori--Bogoliubov information (\cite{AmariNagaoka2000}, Section~7.3) with the Helstrom information, applying to quantum estimation
the large deviation viewpoint of Bahadur \cite{Bahadur1960,Bahadur1967}. Geometric
relations based on the metric properties of quantum information
are studied in \cite{Petz1996,PetzSudar1996,GibiliscoIsola2003,Jencova2004}. Recently, Luo \cite{Luo2006} uncovered the
relation between Helstrom and Wigner--Yanase \cite{WignerYanase1963} skew
information, previously investigated by Luo \cite{Luo2003} and also
considered in
\cite{GibiliscoImparatoIsola2009}, in the most
basic quantum systems, called pure states. Specifically, Luo \cite{Luo2006} showed that, in pure states, Wigner--Yanase information,
$I_{\mathit{WY}}(\theta)$, is exactly twice the Helstrom information,
$I_{\mathit{WY}}(\theta) = 2I_{H}(\theta),$ but left unsolved the problem
of obtaining an exact solution in the generic mixed-state case which
arises when mixtures of pure states are considered.

In this paper, we derive an explicit expression for Wigner--Yanase
information in mixed-state models and relate it to both the
Helstrom information obtained in the mixed state and the
Wigner--Yanase information obtained in the pure states involved in
the mixture. The connection between Wigner--Yanase and Helstrom
information is explored in two-dimensional mixtures of orthogonal
pure states, where a convenient expression for Helstrom
information is available \cite{Luati2004}. While in pure states,
$I_H(\theta)$ and $I_{\mathit{WY}}(\theta)$ are equal up to a constant
factor, in mixed states, their relation depends on the mixing
coefficients. We show that for certain values of the mixing
coefficients, $I_{\mathit{WY}}(\theta)$ is approximately equal to
$I_{H}(\theta)$. This suggests the possibility of approximating
$I_{H}(\theta)$ by $I_{\mathit{WY}}(\theta)$ in the quantum Cram\'{e}r--Rao
bound (\ref{eq:qrb}), as well as in the information inequality
(\ref{eq:infineq}). 
The advantage lies in the fact that
Wigner--Yanase information can be easily calculated, unlike
Helstrom information. In fact, while the latter arises as the
solution of an implicit second order matrix differential equation,
Wigner--Yanase information is simply obtained by scalar
differentiation. Hence, the approximation provides a feasible
bound.

The paper is organized as follows. Section~\ref{sec:cqi} reviews
the basics of quantum statistical inference and introduces the
quantum analogs of classical expected Fisher information.
Further details on probability and statistics in quantum systems
can be found in the books by Holevo \cite{Holevo1982} and Helstrom \cite{Helstrom1976} and
in the paper by Barndorff-Nielsen, Gill and Jupp \cite{BarndorffNielsenGillJupp2003}; excellent
references for quantum information theory are Nielsen and Chuang \cite{NielsenChuang2000} and Petz \cite{Petz2008}. The relation between Helstrom and
Wigner--Yanase information in mixed states is derived in Section
\ref{sec:crb} and the generalizations are discussed in Section
\ref{sec:d}. Proofs and technical details are deferred to the \hyperref[sec:app]{Appendix}. We
conclude this introduction by observing that while quantum
Cram\'{e}r--Rao bounds are intended for vector parameter models,
there is no general relation between the various quantum
information inequalities when the parameter is a vector
\cite{BarndorffNielsenGillJupp2003}. Hence, in the following,
we will restrict our attention to one-parameter models.

\section{Classical and quantum information}\label{sec:cqi}

This paper is concerned with two quantum analogs of classical
expected Fisher information that arise by generalizing to the
quantum setting the following expressions
\begin{eqnarray}
i(\theta,M) & =& \int _{\mathbb{G}_{+}}
\biggl(\frac{\partial}{\partial\theta
} \log p(x;\theta) \biggr)^2p(x;\theta)\mu(\mathrm{d}x)  \label{eq:i1}\\
& =& 4\int _{\mathbb{G}_{+}}
 \biggl(\frac{\partial}{\partial\theta}\sqrt{p(x;\theta)} \biggr)^2
\mu(\mathrm{d}x), \label{eq:i2}
\end{eqnarray}
where $X\dvtx (\Omega,\mathcal{F},P)
\rightarrow(\mathbb{G},\mathcal{G},P_{X})$ is a random variable
characterized by the distribution $P_X(\cdot;\theta)$ with
density $p(x;\theta)$ with respect to the $\sigma$-finite measure
$\mu$ on $(\mathbb{G},\mathcal{G})$, $\mathbb{G} _{+}=\{x\in
\mathbb{G}:p(x;\theta)>0\}$, $\theta\in\Theta\subset\mathbb R$ is
the unknown parameter of interest and the usual regularity
conditions apply \cite{Cramer1946}. The first equation bears a
conceptual meaning, being the expected value of the square score
function $l_{/\theta} = \frac{\partial}{\partial\theta} \log
p(x;\theta)$, while the second is just an equivalent formulation.

If the random variable $X$ describes an experiment in a quantum
system, then its probability law depends on the state of the
system, denoted by $\rho(\theta)$, as it is specified by the
unknown parameter, and on the measurement $M$ that is carried out
in the system, in the following way (trace rule for probability):
\[
P_X(G) = \tr\{\rho(\theta) M(G)\}    \qquad\forall G \in\mathcal
G.
\]
More formally, $\rho(\theta) $ is a density matrix, that is, a
self-adjoint, non-negative and trace-one linear operator acting on
an a $n$-dimensional complex Hilbert space $\mathcal{H}_n$, and $M$
is a probability operator-valued measure, that is, a set of
non-negative self-adjoint linear operators defined on the measure
space $(\mathbb{G},\mathcal{G})$ and taking values in
$\mathcal{H}_{n}$, such that $M(\mathbb{G})=\mathbf{I}$, the
identity operator, $M(\varnothing)=\mathbf{O}$, the null operator,
and
$M(\bigcup _{h=1}^{\infty}G_{h})=\sum _{h=1}^{\infty
}M(G_{h})$ if $G=\bigcup _{h=1}^{\infty}G_{h},G_{h}\cap
G_{l}=\varnothing$, for all $h,l=1,\ldots,\infty$, $h\neq l$. Whenever
$n<\infty$, the Hilbert space $\mathcal H_n$ can be identified
with the $n$-dimensional Euclidean complex space $\mathbb C^n$
endowed with the standard inner product and it is equivalent to
refer to self-adjoint operators or to Hermitian matrices. If $M$ is
absolutely continuous with respect to $\mu$ on
$(\mathbb{G},\mathcal{G})$, such that
$M(G)=\int _{G}m(x)\mu(\mathrm{d}x)$ for all $ G \in\mathcal{G}$,
where $m(x)$ is non-negative and Hermitian, then $P_X(\cdot;
\theta)$ is absolutely continuous with respect to $\mu$ and the
density of $X$ is
%
\begin{equation}p(x;\theta)
=\mbox{tr}\{ \rho(\theta) m(x)\}.\label{eq:p(x,q)}
\end{equation}
The expected value of the random variable $X$ is
\[
\E\{X\}=\int_{\mathbb G} x p(x;\theta)\mu(\mathrm{d}x) =\tr \biggl\{ \rho(\theta)
\int_{\mathbb G}xm(x)\mu(\mathrm{d}x) \biggr\}.
\]
Combining (\ref{eq:i1})
or (\ref{eq:i2}) with (\ref{eq:p(x,q)}), we obtain the classical
expected Fisher information
\[
i(\theta,M) = \int _{\mathbb{G}_{+}} (\mbox{tr}\{\rho
_{/\theta}m(x)\} )^{2} p(x;\theta)^{-1}\mu(\mathrm{d}x),
\]
where $\rho
_{/\theta}$ is the matrix whose $ij$th generic element is the
derivative with respect to $\theta$ of the generic element of
$\rho(\theta)$, that is,
$[\rho_{/\theta}]_{ij}=\frac{\partial}{\partial
\theta}[\rho(\theta)] _{ij}$.

Quantum analogs of classical information are obtained without
performing any measurement on the system. The density matrix
$\rho(\theta)$ plays the role of the density $p(x;\theta)$ and the
expected value of any observable $A$, defined as a self-adjoint
operator, is $\E\{A\} = \tr\{\rho(\theta)A\}$. Note that by defining
$A:=\int_{\mathbb G}xm(x)\mu(\mathrm{d}x)$, we have the formal connection
with the expected value of the random variable $X$.

Helstrom \cite{Helstrom1967} obtained the quantum information $I_H(\theta)$ by
generalizing (\ref{eq:i1}) to the operator-valued quantum setting via
the symmetric logarithmic derivative, which is the
self-adjoint operator $\rho_{/\!\!/\theta}$ implicitly defined by the
relation
%
\begin{equation}
\rho_{/\theta}=\tfrac{1}{2} [\rho(\theta)\rho
_{/\!\!/\theta}+\rho_{/\!\!/\theta}\rho
(\theta) ].\label{eq:sld}
\end{equation}
The result is
%
\begin{equation}
I_H(\theta) = \tr\{\rho(\theta)\rho_{/\!\!/\theta}^{2}\},
\label{eq:IH}
\end{equation}
which, in fact, is the expected value of the
(observable) square symmetric $\rho_{/\!\!/\theta}$. The quantity
$\rho_{/\!\!/\theta}$ is also known as the \textit{quantum score} since, besides
(\ref{eq:IH}), it satisfies
%
\begin{equation}\E\{\rho_{/\!\!/\theta}\}=0 \label{eq:Esld}
\end{equation}
(see \cite{BarndorffNielsenGillJupp2003}, page~789). The proof of the information inequality
(\ref{eq:infineq}) for the one-dimensional parameter case is based
on the Cauchy--Schwarz inequality with Hilbert--Schmidt inner product.
The bound~(\ref{eq:qrb}) can be obtained as a consequence of
(\ref{eq:infineq}) or, directly, by following the derivation of the
classic Cram\'{e}r--Rao bound.

On the other hand, Wigner--Yanase \cite{WignerYanase1963} information, denoted by
$I_{\mathit{WY}}(\theta)$, is obtained by generalizing (\ref{eq:i2}) in a
straightforward manner:
%
\begin{equation}I_{\mathit{WY}}(\theta)=
4\tr \{ [ (\rho(\theta)^{1/2} )_{/\theta
} ]^2 \}.
\label{eq:IWY}
\end{equation}
%

Luo \cite{Luo2006} derived a relation between $I_{\mathit{WY}}(\theta) $ and
$I_{H}(\theta)$ in the case where the system is described by a pure
state density matrix of the form
$\rho(\theta)=|\psi(\theta)\rangle\langle\psi(\theta)|$, where
$|\psi(\theta)\rangle$ is a unit vector in $\mathcal H_n$ and we
have used the Dirac notation according to which
$|\psi(\theta)\rangle$ denotes a column vector and $\langle\psi
(\theta)|$ is its Hermitian transpose. Specifically, the relation~is
%
\begin{equation}I_{WYp}(\theta) = 2
I_{Hp}(\theta),\label{eq:wyhp}
\end{equation}
where the attached
subscript $p$ denotes here that the quantities are associated with a
pure state. The proof can be shortened by observing that, in pure
states, $\rho(\theta)^{1/2}=\rho(\theta)$ and
\[
I_{Hp}(\theta) = 2\tr \{(\rho_{/\theta})^2 \}
\]
(see \cite{Luati2008},~equation (3.1)), from which (\ref{eq:wyhp}) clearly follows by direct
comparison with (\ref{eq:IWY}).

According to quantum state estimation theory, pure states
represent the best knowledge one can have about some specific
properties of the system under observation. Mixed states, obtained
as convex combinations of pure states, indicate a situation of
partial knowledge of the system. They represent probabilistic
mixtures, in the sense that the system under observation is in the
state $%
\rho_{i}(\theta) $ with probability $w_{i}(\theta)$, $i=1,\ldots
,m$, and $\sum_{i=1}^{m}w_i(\theta) = 1$.

To obtain a relation between Helstrom and Wigner--Yanase
information in mixed states, it is necessary to express both of the
information quantities as functions of common variables (up to
some transformation), such as, for example, the corresponding
information quantities in the mixing pure states. In this setting,
the main difficulty is related to Helstrom information since
explicit solutions of (\ref{eq:sld}) when $\rho(\theta)$ is a
mixed state are not usually available in a convenient form.
Restricting to mixtures of two-dimensional orthogonal pure states
allows this difficulty to be overcome and also a
suitable expression for skew information to be derived.

Two-dimensional systems play a crucial role in quantum mechanics.
Electrons, qubits and spin-$\frac{1}{2}$ particles are just some
examples of systems in $\mathbb C^2$. In addition, two-dimensional
mixtures of orthogonal states have an appealing geometric
interpretation, based on the Bloch (or Poincar\'{e}, or Riemann)
sphere representation of states in two-dimensional complex Hilbert
spaces by means of unit vectors in the real
three-dimensional Euclidean space. 
In fact, if $ \mathcal{H}_n=\mathbb{C}^{2}$, then the set of pure
states is the surface of the unit sphere and the set of mixed
states is the interior of the corresponding unit ball (see \cite{Petz2008}). Mixtures of two pure states can be represented as points in
the interior of the sphere, on the straight line joining the two
points on the surface. If the generating pure states are
orthogonal (opposite on the sphere), then the corresponding mixed
states lie on the diameters of the great circles and therefore the
set of such states with given weights can be represented by the
spheres embedded in the unit sphere with the same center, but with radius
less than one and dependent on the weights of the mixtures.
Specifically, we shall consider all of the spheres embedded in the
unit sphere with the same center and radius equal to $|2w(\theta)-1|$,
where $w(\theta)\in(0,1)$ determines the weights of the mixture.

\section{Skew information in mixed states}\label{sec:crb}

A one-parameter two-dimensional mixed state can be represented as
%
\begin{equation} \rho(\theta) =w(\theta)\rho_{1}(\theta)+\bigl(1-w(\theta)\bigr)
\rho_{2}(\theta),\label{eq:ms}
\end{equation}
where $\rho
_{1}(\theta) =| \psi_{1}(\theta)\rangle\langle\psi
_{1}(\theta)| $ and $\rho_{2}(\theta) =| \psi_{2}(\theta
)\rangle\langle\psi_{2}(\theta)| $ are orthogonal pure states,
in the sense that $\langle\psi_1(\theta)|\psi_2(\theta)\rangle=0
$ and $ w(\theta) $ is a function of $\theta$ taking values in
the real interval $(0,1)$, $w(\theta)\neq\frac{1}{2}$. Note that
orthogonal pure states are such that
$\rho_h(\theta)\rho_h(\theta)=\rho_h(\theta)$ and
$\rho_h(\theta)\rho_k(\theta) =\mathbf O$ for $h\neq k$. Let us
denote by $\rho_{/\theta h }$ and $\rho_{/\!\!/\theta h }$ the
term-by-term first derivative and the symmetric logarithmic
derivative of $\rho_{h}(\theta), h=1,2$, with respect to
$\theta$, respectively, and by $I_{Hh}(\theta)$ and
$I_{WYh}(\theta)$, the Helstrom and Wigner--Yanase information
extracted from the pure state $\rho_h(\theta)$, respectively. We assume that
$w(\theta)$ has continuous first derivative, $w_{/\theta}$, such
that $w_{/\theta}\rightarrow0$ faster than $\sqrt{w(\theta)}$ and
$\sqrt{1-w(\theta)}$ for all $\theta$ such that
$w(\theta)\rightarrow0$ and $w(\theta)\rightarrow1$,
respectively.
With these premises, we can state the following proposition.
\begin{proposition}\label{prop1}
In the mixed state
(\ref{eq:ms}), Wigner--Yanase skew information is
%
\begin{equation}
I_{\mathit{WY}}(\theta) = \frac{( w_{/\theta}) ^{2}}{w( \theta) (
1-w( \theta) ) }+ \bigl(1-
2\sqrt{w(\theta)\bigl(1-w(\theta)\bigr)} \bigr)I_{\mathit{WY}1}( \theta
)\label{eq:IWYms}
\end{equation}
and the following relation holds
%
\begin{equation}I_{\mathit{WY}}(\theta) =
\alpha_w (\theta) I_H(\theta) + \beta_w (\theta),
\label{eq:IWYandIH}
\end{equation}
where
\[
\alpha_w (\theta)=
\frac{2}{1+2\sqrt{w(\theta)(1-w(\theta))}} \quad\mbox{and}\quad \beta_w
(\theta)= -
\frac{1-2\sqrt{w(\theta)(1-w(\theta))}}{1+2\sqrt{w(\theta)(1-w(\theta))}}
\frac{(w_{/\theta})^2}{w(\theta)(1-w(\theta))}.
\]
\end{proposition}

The proof is in Section \ref{sec:app1} and, concerning
(\ref{eq:IWYms}), it is based on some properties of pure states
and quantum information quantities. The proof of
(\ref{eq:IWYandIH}) is based on the comparison with Helstrom
information derived in \cite{Luati2004}, Lemma~2, where a specific
choice of $|\psi_1(\theta)\rangle$ and $|\psi_2(\theta)\rangle$
was made in order to derive the symmetric logarithmic derivative
of the mixed state as a function of the symmetric logarithmic
derivative of $\rho_1(\theta)$,
%
\begin{equation}I_H(\theta) =\frac{( w_{/\theta}) ^{2}}{w( \theta) ( 1-w(
\theta) ) }+\bigl( 2w( \theta) -1\bigr) ^{2}I_{H1}( \theta
).\label{eq:IHms}
\end{equation}
The boundary conditions on
$w_{/\theta}$ ensure that $I_{\mathit{WY}}(\theta) \rightarrow
I_{\mathit{WY}1}(\theta)$ and $I_{H}(\theta) \rightarrow I_{H1}(\theta)$
when $w(\theta)\rightarrow0,1$ and imply that $\beta_w
(\theta)\rightarrow0$ when $w (\theta)\rightarrow0,1$. Hence,
they encompass the case where the mixing coefficient does not
depend on $\theta$ ($w_{/\theta}=0$ $\forall\theta\in\Theta\subset\mathbb R$).

In the analysis of the relation between $I_{\mathit{WY}}(\theta)$ and
$I_{H}(\theta)$, we first consider the latter case, that is, $w$ does
not depend on $\theta$, which implies
that
\[
I_{\mathit{WY}}(\theta) = \alpha_w I_H(\theta),
\]
where $\alpha_w$ is
function of $w$, symmetric with respect to $w=\frac{1}{2}$, where
it reaches its minimum of~1, and with maximum of 2 at
$w=0$ and $w=1$. The implications are evident: when the mixing
coefficients tend to $0$ and $1$, the skew information tends
to be twice the Helstrom information since the mixture tends to
reproduce a pure state system. Conversely, when
$w\rightarrow\frac{1}{2}$, the distance between the two quantum
information quantities reaches its minimum, in that they tend to
be equal. The limit case $w\rightarrow\frac{1}{2}$ is of special
interest to us, so we shall return to it in the discussion of
Proposition~\ref{prop1} for the case where $w$ depends on $\theta$.
Furthermore, both $I_{\mathit{WY}}(\theta)$ and $I_H(\theta)$ are smaller
than the corresponding information quantities in the pure states.
This is coherent with quantum theory, according to which pure
states represent the best knowledge that one can have about a
quantum system.

On the other hand, if the mixing coefficients depend on $\theta$,
then a sufficient condition for $I_{\mathit{WY}}(\theta)<I_{\mathit{WY}1}(\theta)$
and $I_{H}(\theta)<I_{H1}(\theta)$ is that
\[
I_{\mathit{WY}1}(\theta)>\frac{( w_{/\theta}) ^{2}}{2 w^2( \theta) (
1-w( \theta) )^2},
\]
which follows by simple algebraic
manipulation of (\ref{eq:IWYms}) and (\ref{eq:IHms}), together with the identity
$I_{\mathit{WY}1}(\theta) = 2I_{H1}(\theta)$. In other words, it is
sufficient that the condition $I_{H}(\theta)<I_{H1}(\theta)$ is
satisfied to ensure that $I_{\mathit{WY}}(\theta)<I_{\mathit{WY}1}(\theta)$.

The limit cases $w(\theta)\rightarrow0,\frac{1}{2},1$ turn out to
be identical to the limit cases when $w$ is constant. In
particular, if $w(\theta)\rightarrow\frac{1}{2}$, then
$\beta_w(\theta)\rightarrow0$ and $\alpha_w(\theta)\rightarrow
1$ so that we find $I_{\mathit{WY}}(\theta)\rightarrow I_H(\theta)$.
Hence, we can state the following corollary to Proposition~\ref{prop1}.
This is a direct consequence of the quantum Cram\'{e}r--Rao bound
(\ref{eq:qrb}) and hence the proof is omitted.
\begin{corollary}\label{cor1}
In the mixed state
(\ref{eq:ms}), if $w(\theta)\rightarrow\frac{1}{2}$, then
$I_{\mathit{WY}}(\theta)\rightarrow I_{H}(\theta)$ and
$I_{\mathit{WY}}(\theta)^{-1}$ constitutes an approximate lower bound for
$\Var\{t(x)\}$.
\end{corollary}

The equality between $I_{H}(\theta)$ and $I_{\mathit{WY}}(\theta)$ is a
limit condition that holds in a degenerate case. Precisely, if
$w(\theta) = \frac{1}{2}$, then $\rho(\theta)=
\frac{1}{2}\mathbf{I}$. This case (the center of the unit sphere)
represents the maximum entropy situation, that is, complete ignorance
about the quantum system under observation. Nevertheless, the
bound has the following interpretation: the precision of the
feasible Wigner--Yanase quantum Cram\'{e}r--Rao bound increases as
long as the mixture approaches the maximum entropy case. Values of
$w(\theta)$ that lie in a neighborhood of $w(\theta)=\frac{1}{2}$
give rise to sensible mixtures and feasible approximated quantum
Cram\'{e}r--Rao bounds. In the two-dimensional case, the size of
the approximation can be measured using equation
(\ref{eq:IWYandIH}).

In pure states, where we are equipped with all the quantum
information given by the system, knowledge of $I_{\mathit{WY}}$ is equivalent to
knowledge of $I_H$ (up to a constant factor, they are equal); in
mixed states, the distance between the two quantum information
quantities depends on the probability of being in a
(one-dimensional) space rather than on its orthogonal complement.
The same distance between $I_H(\theta)$ and $I_{\mathit{WY}}(\theta)$
decreases as long as the knowledge of the system diminishes,
eventually tending to zero in the maximum entropy case (uniform
distribution on the coefficients).

\section{Discussion }\label{sec:d}

We derived the relation between Helstrom and Wigner--Yanase
information in two-dimensional mixed-state systems, where a geometric
interpretation of the states as spheres embedded in the unit
sphere can be drawn. We now discuss the case of general mixed
states
\[
\rho(\theta) = \sum_{i=1}^{m}w_i(\theta)\rho_i(\theta),
\]
where
the $\rho_i(\theta)$ are pure states defined on $\mathbb C^n$ and
not necessarily orthogonal, and
\[
\sum_{i=1}^m w_i(\theta) = 1.
\]
If the mixed state admits the spectral decomposition
%
\begin{equation}\rho(\theta) = \sum_{l=1}^n
\lambda_l(\theta)\rho_l (\theta),\label{eq:sd}
\end{equation}
then
$I_{\mathit{WY}}(\theta)$ can be obtained directly, following the proof of
Proposition~\ref{prop1}, as
%
\begin{equation}I_{\mathit{WY}}(\theta) = \sum_{l=1}^n \lambda_l(\theta)
I_{WY,l}(\theta) + \sum_{l=1}^n \frac{(\lambda_{/\theta
l})^2}{\lambda_l(\theta)} + 4 \sum_{l=1}^n \sum_{k\neq l}
\sqrt{\lambda_l(\theta)\lambda_k(\theta)}\tr\{\rho_{/\theta
l}\rho_{/\theta k}\}.\label{eq:IWYmsgen}
\end{equation}
When $n=2$, $\lambda_1(\theta) = w(\theta)$, $\lambda_2(\theta) =
1-w(\theta)$, $(\lambda_{/\theta1})^2=(\lambda_{/\theta2})^2
=(w_{/\theta})^2$ and (\ref{eq:IWYmsgen}) is equal to
(\ref{eq:IWYms}). On the other hand, concerning $I_{H}(\theta)$, a
solution of (\ref{eq:sld}) based on the spectral decomposition of
$\rho(\theta)$ is \cite{PetzSudar1996}
%
\begin{equation}\rho_{/\!\!/\theta} = \sum_{l=1}^n
\sum_{k=1}^n\frac{2}{\lambda_l(\theta)+\lambda_k(\theta)}\rho_j(\theta
)\rho_{/\theta}\rho_k(\theta),
\label{eq:sld_sd}
\end{equation}
which allows us to obtain the
maximum Fisher information attainable in the mixed state
(\ref{eq:sd}),
%
\begin{equation}I_H(\theta) = \sum_{l=1}^n \frac{(\lambda_{/\theta
l})^2}{\lambda_l(\theta)} + \sum_{l=1}^n \sum_{k\neq l}
\sum_{z=1}^n \frac{4 \lambda_l(\theta)
(\lambda_k(\theta)-\lambda_l(\theta) )\lambda_z (\theta)}
{(\lambda_l(\theta)+\lambda_k(\theta))^2} \tr \{\rho_l(\theta)\rho
_{/\theta k}\rho_{/
\theta z} \}.\label{eq:IHmsgen}
\end{equation}
The derivation is not direct and is therefore deferred to
Section \ref{sec:app4} of the \hyperref[sec:app]{Appendix}. It is easy to verify that when $n=2$,
(\ref{eq:IHmsgen}) is equal to (\ref{eq:IHms}), since, besides
$\lambda_1(\theta) = w(\theta)$, $\lambda_2(\theta) = 1-w(\theta)$
and $(\lambda_{/\theta1})^2=(\lambda_{/\theta2})^2
=(w_{/\theta})^2$, $4 \tr\{\rho_l(\theta)\rho_{/\theta k}\rho_{/
\theta z}\}$ is equal to plus or minus
$I_{H1}(\theta)$, according to whether $z=k$ or $ z \neq k$, respectively,
as follows by the properties of two orthogonal pure states in
two-dimensional spaces, $\rho_{/\!\!/\theta h} = 2 \rho_{/\theta h}$,
$h=1,2$, $\rho_{/\theta h} = -\rho_{/\theta k} $, $h \neq k$ and
$I_{H2}(\theta) = I_{H1}(\theta)$ (\cite{Luati2004}, equations (i), (vi) and
(ix), respectively).

It is evident that the two information quantities
(\ref{eq:IWYmsgen}) and (\ref{eq:IHmsgen}) are sensibly comparable
only under some specific assumptions on the state of the quantum
system, like those we have made in Section \ref{sec:crb}.
Moreover, the two expressions do not explicitly depend on the
mixing coefficients, nor on the mixing pure states, so their
formal equivalence loses its interpretation. Nevertheless,
rewriting $I_{\mathit{WY}}(\theta)$ in (\ref{eq:IWYmsgen}) as a function
of $I_H(\theta)$ in (\ref{eq:IHmsgen}) provides an
eigenvalue-based condition of equality between the two quantum
information quantities, one which allows us to derive a result
equivalent to Corollary~\ref{cor1} in Section \ref{sec:crb}, that is, an
approximate quantum Cram\'{e}r--Rao bound based on skew
information in $n$-dimensional mixed states that admit the
spectral decomposition (\ref{eq:sd}). We state this finding in a
proposition, proved in Section \ref{sec:app3} of the \hyperref[sec:app]{Appendix}, and a corollary,
whose proof is omitted since it follows from the proposition in a
straightforward manner.
\begin{proposition}\label{prop2}
In the mixed state
(\ref{eq:sd}),
\begin{eqnarray*}
I_{\mathit{WY}}(\theta) = I_{H}(\theta) + \gamma_{\lambda,\rho
}(\theta),
\end{eqnarray*}
where, omitting the arguments of the functions,
%
\begin{equation}\gamma_{\lambda,\rho} = - 4\sum_{l=1}^n\sum_{k\neq l}
 \Biggl(
\bigl(\lambda_l-\sqrt{\lambda_l\lambda_k}\bigr)\tr\{\rho_{/l}\rho_{/ k}\} +
\sum_{z=1}^n \frac{\lambda_l (\lambda_k-\lambda_l )\lambda_z }
{(\lambda_l+\lambda_k)^2} \tr \{\rho_l\rho_{/ k}\rho_{/
z} \} \Biggr).\label{eq:gamma}
\end{equation}
\end{proposition}

The proof essentially consists of rewriting Wigner--Yanase
information in a convenient way and then comparing it with
Helstrom information. The quantity $\gamma_{\lambda,\rho}(\theta)$
depends both on the eigenvalues of $\rho(\theta)$, playing the
role of the mixing coefficients in the representation
(\ref{eq:sd}), and on the states $\rho_l(\theta)$, via their
derivatives. Hence, it would be misleading to interpret
$\gamma_{\lambda,\rho}(\theta)$ as a pure additive quantity, be it
positive or negative. As a matter of fact, the dependence of
$\gamma_{\lambda,\rho}(\theta)$ on $\lambda_l(\theta)$ and
$\rho_l(\theta)$ implies that the quantity
$\gamma_{\lambda,\rho}(\theta)$ can involve, at least,
combinations of information quantities in the one-dimensional pure
states that span the mixed-state system (\ref{eq:sd}). To illustrate
this, we shall consider the two-dimensional case, where
$\gamma_{\lambda,\rho}(\theta) =
 (1-2\sqrt{w(\theta)(1-w(\theta))} )^2 I_{H1} (\theta) =
(\alpha_w(\theta)-1) I_H(\theta) +\beta_w(\theta)$. Unlike $\gamma
_{\lambda,\rho}(\theta)$, the quantities
$\alpha_w(\theta)$ and $\beta_w(\theta)$ depend only on the mixing
coefficients and their interpretation is unequivocal. Also, note
that in the two-dimensional case, $\gamma_{\lambda,\rho}(\theta)
\rightarrow0$ for all $\theta$ such that $w(\theta) \rightarrow
\frac{1}{2}$, this leading to the approximate quantum
Cram\'{e}r--Rao bound stated in Corollary~\ref{cor1}.

Corollary~\ref{cor1} can be generalized to the $n$-dimensional case by
means of Proposition~\ref{prop2}. In fact, it is immediate to see that in
(\ref{eq:gamma}), if $\lambda_k(\theta)\rightarrow
\lambda_l(\theta)$ for all $k$ and $l$, then
$\gamma_{\lambda,\rho}(\theta) \rightarrow0$ and
$I_{WY(\theta)}\rightarrow I_H(\theta)$. A sufficient condition
for all of the $\lambda$ to be equal is that
$\lambda_l(\theta)\rightarrow\frac{1}{n}$ for all $l=1,\ldots,n$
and the following corollary can be established as a consequence of
the quantum Cram\'{e}r--Rao bound~(\ref{eq:qrb}).
\begin{corollary}\label{cor2}
In the mixed state
(\ref{eq:sd}), if $\lambda_l(\theta)\rightarrow\frac{1}{n}$
for all $l=1,2,\ldots, n$, then $I_{\mathit{WY}}(\theta)\rightarrow
I_{H}(\theta)$ and $I_{\mathit{WY}}(\theta)^{-1}$ constitutes an
approximate lower bound for $\Var\{t(x)\}$.
\end{corollary}

We can therefore conclude that in a mixed state specified by its
spectral decomposition, a limit condition analog to the one
holding in two-dimensional systems holds as well. The
interpretation of the approximate quantum Cram\'{e}r--Rao bound
specified in Corollary~\ref{cor2} is the same as in the two-dimensional
case. As long as, in (\ref{eq:sd}), the mixing coefficients tend
to be uniformly distributed over $n$, Wigner--Yanase skew
information tends to be equal to Helstrom information and can
serve as an approximation of the latter in the information
inequality (\ref{eq:infineq}) and in the quantum Cram\'{e}r--Rao
bound (\ref{eq:qrb}). The uniform distribution of the coefficients
is a limit condition that holds in the degenerate maximum entropy
case, but situations where the limit is approached give rise to
sensible mixtures and feasible approximated quantum Cram\'{e}r--Rao
bounds. Deriving equations~(\ref{eq:IWYmsgen}) and
(\ref{eq:IHmsgen}) shows that it is, in general, much easier to
obtain a feasible expression for Wigner--Yanase information
than for Helstrom information, which further justifies the use of~(\ref{eq:IWY}) as an approximate upper bound for the classical
expected Fisher information.

\begin{appendix}
\section*{Appendix}\label{sec:app}
\renewcommand{\theequation}{\arabic{equation}}
\subsection{\texorpdfstring{Proof of Proposition~\protect\ref{prop1}}{Proof of Proposition 1}}\label{sec:app1}
We derive $I_{\mathit{WY}}(\theta) =
4\tr\{[(\rho(\theta)^{1/2})_{/\theta})]^2\}$ in the mixed
state (\ref{eq:ms}). Since $|\psi_1(\theta)\rangle$ and
$|\psi_2(\theta)\rangle$ are orthogonal and $\rho_1(\theta)$ and
$\rho_2(\theta)$ are pure states,
we have
\begin{equation} \rho(\theta)^{1/2}
= \sqrt{w(\theta)} \rho_1(\theta) + \sqrt{1-w(\theta)}
\rho_2(\theta). \label{eq:sqrt}
\end{equation}
Note that
$\rho(\theta)^{1/2}\rightarrow\rho_h(\theta)$, $h=2,1$,
if $w(\theta)\rightarrow0,1$. Also, note that
$\rho(\theta)^{1/2}$ is not a density matrix since its
trace is not equal to 1. The choice of
$\rho(\theta)^{1/2}$, as in (\ref{eq:sqrt}), as opposed its
alternative with the negative sign, guarantees that the square
root of the positive semidefinite self-adjoint operator
$\rho(\theta)$ is positive semidefinite for all of the values
$w(\theta)$ (\cite{HornJohnson1991}, Theorem 2.6, page 405).

The elementwise derivative of (\ref{eq:sqrt}) with respect to
$\theta$ is given by
\begin{eqnarray*} (\rho(\theta)^{1/2} )_{/\theta} =
\frac{w_{/\theta}}{2\sqrt{w(\theta)}} \rho_1(\theta) +
\sqrt{w(\theta)} \rho_{/\theta1} -
\frac{w_{/\theta}}{2\sqrt{1-w(\theta)}} \rho_2(\theta) +
\sqrt{1-w(\theta)} \rho_{/\theta2},
\end{eqnarray*}
where we have assumed that $w_{/\theta}$ tends
to zero faster than $\sqrt{w(\theta)}$ and $\sqrt{1-w(\theta)}$
for all $\theta$ such that $w(\theta)\rightarrow0,1$,
respectively.

In taking the square and then the trace of the above quantity, we
observe that the spectral theorem in $\mathbb C^2$ implies that
$\rho_2(\theta) = \bI-\rho_1(\theta)$ and, consequently, that
$\rho_{/\theta2} = - \rho_{/\theta1}$, $\rho_{/\theta
2}\rho_{/\theta1} = -(\rho_{/\theta1})^2$ and $(\rho_{/\theta
2})^2 = (\rho_{/\theta1})^2 =
((\rho_1(\theta)^{1/2})_{/\theta})^2$. Furthermore,
%
\begin{equation}\tr\{\rho_k(\theta)\rho_{/\theta h}\}=0 \qquad  \forall
h,k.\label{eq:tr_r_r/}
\end{equation}
In fact, when $h=k$, equation (\ref{eq:tr_r_r/})
follows from the fact that, in pure states, $\rho_{/\!\!/\theta} =
2\rho_{/\theta}$ (\cite{Luati2008}, equation (3.1)), combined with the fact that
$\E\{\rho_{/\!\!/\theta}\}= 2\E\{\rho_{/\theta}\} = 0$, stated in
(\ref{eq:Esld}). On the other hand, when $h\neq k$, using the
definition of symmetric logarithmic derivative (\ref{eq:sld}), $
\tr\{\rho_k(\theta) \rho_{/\theta h}\} =
\frac{1}{2}\tr\{(\rho_k(\theta)\rho_h(\theta)\rho_{/\!\!/\theta
h}+\rho_k(\theta)\rho_{/\!\!/\theta h}\rho_h(\theta))\}=0$.

Hence,
\begin{eqnarray*}
4\tr \{ [ (\rho(\theta)^{1/2})_{/\theta} ]^2 \}
&= &   4 \biggl[ \biggl(\frac{w_{/\theta}}{2\sqrt{w(\theta)}} \biggr)^2
+ \biggl(\frac{w_{/\theta}}{2\sqrt{1-w(\theta)}}  \biggr)^2 \biggr]\\
 &&{} + 4 \bigl(w(\theta) + 1 - w(\theta) -
2\sqrt{w(\theta)}\sqrt{1-w(\theta)} \bigr) \tr\{(\rho_{/\theta1})^2\}
\\
&= & \frac{( w_{/\theta}) ^{2}}{w( \theta) ( 1-w( \theta) )
}+ \bigl(1- 2\sqrt{w(\theta)\bigl(1-w(\theta)\bigr)} \bigr)I_{\mathit{WY}1}( \theta).
\end{eqnarray*}
This completes the first part of the proof, that is, equation
(\ref{eq:IWYms}).

We now prove equation (\ref{eq:IWYandIH}). Without loss of
generality, we can choose
\[
|\psi_2(\theta) \rangle= 2\sqrt{2}
I^{-1/2}_{\mathit{WY}1}(\theta)\rho_{/\theta
1}|\psi_1(\theta)\rangle=
I^{-1/2}_{H1}(\theta)\rho_{\//\theta
1}|\psi_1(\theta)\rangle.
\]
This setting is convenient because
Lemma 2 in \cite{Luati2004} applies, which states that in
two-dimensional orthogonal mixed states where
$|\psi_1(\theta)\rangle$ and $|\psi_2(\theta)\rangle$ are defined
as above, Helstrom information is given by equation
(\ref{eq:IHms}), which can be written as
\[
\bigl( 2w( \theta) -1\bigr) ^{-2}  \biggl[ I_H(\theta) - \frac{( w_{/\theta
}) ^{2}}{w( \theta) ( 1-w( \theta) ) } \biggr] =I_{H1}( \theta
).
\]
It follows from (\ref{eq:wyhp}) and (\ref{eq:IWYms}) that
\[
I_{\mathit{WY}}(\theta) = \frac{(
w_{/\theta}) ^{2}}{w( \theta) ( 1-w( \theta) ) }+ 2 \bigl(1-2\sqrt{w(\theta)\bigl(1-w(\theta)\bigr)}\bigr)I_{H1}( \theta).
\]
Combining the
two former equations, we get
\begin{eqnarray*}
I_{\mathit{WY}}(\theta) &= &\frac{( w_{/\theta}) ^{2}}{w( \theta) ( 1-w( \theta) ) }
 \bigg[1-\frac{ 2 (1- 2\sqrt{w(\theta)(1-w(\theta))} }{( 2w( \theta
) -1) ^{2}} \bigg]
\\
&&{}+ \frac{ 2 (1- 2\sqrt{w(\theta)(1-w(\theta))}
}{( 2w( \theta) -1) ^{2}} I_H(\theta)  \\
&= & \frac{( w_{/\theta
}) ^{2}}{w( \theta) ( 1-w( \theta) ) }  \bigg[\frac{ - ( 1-2
\sqrt{w(\theta)(1-w(\theta))} )^2}{ (2w( \theta)
-1)^{2}} \bigg]
\\
&&{}+ \frac{ 2  (1- 2\sqrt{w(\theta)(1-w(\theta))}
 ) }{( 2w( \theta) -1) ^{2}} I_H(\theta).
\end{eqnarray*}
Given
that $( 2w( \theta) -1) ^{2} =  (1-
2\sqrt{w(\theta)(1-w(\theta))}  )  (1+
2\sqrt{w(\theta)(1-w(\theta))}  ) $, we obtain
\begin{eqnarray*}
I_{\mathit{WY}}(\theta) &=& -\frac{( w_{/\theta})
^{2}}{w( \theta) ( 1-w( \theta) ) } \frac{ 1- 2
\sqrt{w(\theta)(1-w(\theta))}} {1 + 2\sqrt{w(\theta)(1-w(\theta))}
}
\\
&&{}+ \frac{ 2 }{1+ 2\sqrt{w(\theta)(1-w(\theta))}} I_H(\theta).
\end{eqnarray*}

\subsection{\texorpdfstring{Proof of equation (\protect\ref{eq:IHmsgen})}{Proof of equation (19)}}\label{sec:app4}

It follows from (\ref{eq:sd}) and (\ref{eq:sld_sd}) that
\begin{eqnarray*}
I_H(\theta) &=& \tr\{\rho(\theta)\rho^2_{/\!\!/\theta} \}
\\
&=&\tr \Biggl\{
\sum_{l=1}^n \sum_{k=1}^n\sum_{j=1}^n \frac{4\lambda_l(\theta)}
{(\lambda_l(\theta)+\lambda_k(\theta))(\lambda_k(\theta)+\lambda
_j(\theta))}
\rho_l(\theta)\rho_{/\theta}\rho_k(\theta)\rho_{/\theta}\rho_j(\theta
) \Biggr\}.
\end{eqnarray*}
The linearity and the cyclical property of the trace operator,
associated with the orthogonality of pure states $\rho_l(\theta)$
(by the spectral theorem, all of these operators are orthogonal
projections onto one-dimensional supspaces of $\mathbb C^n$), give
\[
I_H(\theta) =
\sum_{l=1}^n \sum_{k=1}^n \frac{4\lambda_l(\theta)}
{(\lambda_l(\theta)+\lambda_k(\theta))^2}\tr \{
\rho_l(\theta)\rho_{/\theta}\rho_k(\theta)\rho_{/\theta} \}
\]
and, omitting all of the arguments,
\begin{eqnarray*}I_H(\theta) &= &
\sum_{l=1}^n \sum_{k=1}^n \frac{4\lambda_l}
{(\lambda_l+\lambda_k)^2}\tr \Biggl\{
\rho_l\sum_{m=1}^n (\rho_m\lambda_{ / m}+\rho_{/ m}\lambda_m)
\rho_k\sum_{z=1}^n (\rho_z\lambda_{ / z}+\rho_{/
z}\lambda_z) \Biggr\} \\
&= & \sum_{l=1}^n \sum_{k=1}^n \frac{4\lambda_l}
{(\lambda_l+\lambda_k)^2}\tr \Biggl\{ \Biggl(
\rho_l\lambda_{/l}+\sum_{m=1}^n \lambda_m \rho_l\rho_{/ m} \Biggr)
 \Biggl( \rho_k\lambda_{/k}+\sum_{z=1}^n \lambda_z \rho_k\rho_{/
z} \Biggr) \Biggr\} \\
&= & \sum_{l=1}^n \frac{(\lambda_{/l})^2}
{\lambda_l} + \sum_{l=1}^n \frac{\lambda_{/l}}
{\lambda_l} \sum_{z=1}^n \lambda_z \tr\{\rho_l\rho_{/
z}\} + \sum_{l=1}^n \frac{\lambda_{/l}}
{\lambda_l} \sum_{m=1}^n \lambda_m \tr\{\rho_l\rho_{/
m}\}  \\
&&{} + \sum_{l=1}^n \sum_{k=1}^n \frac{4\lambda_l}
{(\lambda_l+\lambda_k)^2} \tr \Biggl\{ \sum_{m=1}^n\sum_{z=1}^n \lambda
_m \lambda_z \rho_l\rho_{/
m}\rho_k\rho_{/ z} \Biggr\}  \\
&= & \sum_{l=1}^n \frac{(\lambda_{/l})^2}
{\lambda_l} + 2\sum_{l=1}^n \frac{\lambda_{/l}}
{\lambda_l} \sum_{z=1}^n \lambda_z \tr\{\rho_l\rho_{/
z}\}
\\
&&{}+ \sum_{l=1}^n \sum_{k=1}^n \frac{4\lambda_l}
{(\lambda_l+\lambda_k)^2} \tr \Biggl\{ \sum_{m=1}^n\sum_{z=1}^n \lambda
_m \lambda_z \rho_l\rho_{/
m}\rho_k\rho_{/ z} \Biggr\}.
\end{eqnarray*}

The second term of the above equation is null; see
(\ref{eq:tr_r_r/}). In order to simplify the final expression of
$I_H(\theta)$, it is convenient to analyze the implications of the
product $\rho_l\rho_{/m}$ for the third term of the above equation.
When $l=m$, we can prove that
%
\begin{equation}\tr\{\rho_l\rho_{/l}\rho_l\rho_{/l}\}=0  \qquad\mbox{for
all }
l=1,\ldots,n.\label{eq:tr_null}
\end{equation}
By definition of
pure state, $\rho_l\rho_l = \rho_l$ and, consequently,
$\rho_{/l}\rho_l+\rho_l\rho_{/l}=\rho_{/l}$, that is, $
\rho_l\rho_{/l}= \rho_{/l}-\rho_{/l}\rho_l$, which, substituted into
$\tr\{\rho_l\rho_{/ l}\rho_l\rho_{/ l}\}$, gives
$\tr\{(\rho_{/l}-\rho_{/l}\rho_l)\rho_l\rho_{/l}\} =0$. On the
other hand, let us consider any $m \neq l$. It is straightforward
to prove that
\begin{equation}\rho_l\rho_{/m} = -\rho_{/l}\rho_m,
\qquad m\neq l. \label{eq:rholrhom}
\end{equation}
In fact, from the
spectral theorem, $\rho_l\rho_m =\mathbf{O} \Rightarrow
\rho_{/l}\rho_m + \rho_l\rho_{/m}= \mathbf{O} \Rightarrow
\rho_l\rho_{/m} = -\rho_{/l}\rho_m $, that is, non-null terms arise
for $m\neq l$ when $m=k=z$.

Let us now consider the third term of the last expression for
$I_H(\theta)$. First, note that if $k=l$, then the trace is null.
In fact, if $m=l$, then $\tr \{\rho_l\rho_{/ m}\rho_l\rho_{/
z} \} = \tr \{\rho_l\rho_{/ l}\rho_l\rho_{/ z} \} =
0 $ both if $z=l$, by (\ref{eq:tr_null}), and if $z \neq l$, by
(\ref{eq:rholrhom}), and by the cyclical properties of the trace
operator; on the other hand, if $m\neq l$, then
$\tr \{\rho_l\rho_{/ m}\rho_l\rho_{/ z} \} =
\tr \{-\rho_{/l}\rho_{m}\rho_l\rho_{/ z} \} = 0$ for all~$z$.

Therefore, it follows that $I_H(\theta)$ can be simplified as
\begin{eqnarray*}I_H(\theta) &= &
\sum_{l=1}^n \frac{(\lambda_{/l})^2}
{\lambda_l} + \sum_{l=1}^n \sum_{k\neq l} \frac{4\lambda_l}
{(\lambda_l+\lambda_k)^2} \tr \Biggl\{ \sum_{m=1}^n\sum_{z=1}^n \lambda
_m \lambda_z \rho_l\rho_{/
m}\rho_k\rho_{/ z} \Biggr\}  \\
&= & \sum_{l=1}^n \frac{(\lambda_{/l})^2}
{\lambda_l} + \sum_{l=1}^n \sum_{k\neq l} \frac{4\lambda_l}
{(\lambda_l+\lambda_k)^2} \tr \Biggl\{ \sum_{m\neq k}\sum_{z=1}^n \lambda
_m \lambda_z \rho_l\rho_{/
m}\rho_k\rho_{/ z} \Biggr\} \\
&&{} + \sum_{l=1}^n \sum_{k\neq l} \frac{4\lambda_l}
{(\lambda_l+\lambda_k)^2} \tr \Biggl\{ \sum_{z=1}^n \lambda_k
\lambda_z \rho_l\rho_{/ k}\rho_k\rho_{/ z} \Biggr\},
\end{eqnarray*}
where we have distinguished the case $m\neq k$ from the case $m =
k$. Now, note that, when $m\neq k$, the trace in the second summand
is null unless $m=l$, as follows by (\ref{eq:rholrhom}). Using the
same arguments, but applied to $l\neq k$ in the third summand, we
can simplify the final expression of Helstrom information in a
generic mixed state as

\begin{eqnarray*} I_H(\theta) &= & \sum_{l=1}^n
\frac{(\lambda_{/l})^2}
{\lambda_l} + \sum_{l=1}^n \sum_{k\neq l} \frac{4 (\lambda_l)^2}
{(\lambda_l+\lambda_k)^2} \tr \Biggl\{ \sum_{z=1}^n \lambda_z (- \rho
_l\rho_{/k}\rho_{/
z}) \Biggr\} \\
&&{} + \sum_{l=1}^n \sum_{k\neq l} \frac{4\lambda_l \lambda_k}
{(\lambda_l+\lambda_k)^2} \tr \Biggl\{ \sum_{z=1}^n \lambda_z
(-\rho_{/l}\rho_k\rho_{/z}) \Biggr\}  \\
&= & \sum_{l=1}^n \frac{(\lambda_{/l})^2}
{\lambda_l} + \sum_{l=1}^n \sum_{k\neq l} \sum_{z=1}^n \frac{4 \lambda
_l (\lambda_k-\lambda_l )\lambda_z }
{(\lambda_l+\lambda_k)^2} \tr \{\rho_l\rho_{/k}\rho_{/
z} \}.
\end{eqnarray*}
%

\subsection{\texorpdfstring{Proof of Proposition~\protect\ref{prop2}}{Proof of Proposition 2}}\label{sec:app3}

Let us consider equation (\ref{eq:IWYmsgen}). We observe that
$I_{WY,l}(\theta) = 4\tr\{\rho_{/\theta l}\rho_{/\theta l}\} =
-4\sum_{k\neq l}\tr\{\rho_{/\theta l}\rho_{/\theta k}\}$, as
follows by $\rho_{/\theta l} = -\sum_{k\neq l}\rho_{/\theta k}$, a
consequence of the spectral identity $\sum_{l=1}^n
\rho_l(\theta)=\bI$. Hence,
\[
I_{\mathit{WY}}(\theta) = \sum_{l=1}^n
\frac{(\lambda_{/\theta
l})^2}{\lambda_l(\theta)}-4\sum_{l=1}^n\sum_{k\neq l}
 \bigl(\lambda_l(\theta)-\sqrt{\lambda_l(\theta)\lambda_k(\theta)} \bigr)
\tr\{\rho_{/\theta k}\rho_{/\theta l}\}.
\]
By replacing, in the
above equation, the first term on the right-hand side by the same
quantity obtained from (\ref{eq:IHmsgen}) as a function of
$I_H(\theta)$, and rearranging, we obtain the claimed result.
\end{appendix}

\section*{Acknowledgements}
I would like to thank the Editor and a referee for their comments.

\printhistory

\end{document}